\documentclass[10pt,twoside]{article}
\usepackage{amssymb}
\usepackage{latexsym}
\newtheorem{theorem}{Theorem}[section]

\newtheorem{lemma}{Lemma}[section]

\newtheorem{remark}{Remark}[section]

\def\proof{\mbox {\it Proof.~}}
\makeatletter\def\theequation{\arabic{section}.\arabic{equation}}\makeatother

\usepackage{amsfonts,amsmath}
\usepackage{a4,amssymb}
\usepackage[dvips]{graphicx}
\usepackage[latin1]{inputenc}

\usepackage{color}
\ifx\pdfoutput\undefined 
\else
\fi

\def\Xint#1{\mathchoice
{\XXint\displaystyle\textstyle{#1}}%
{\XXint\textstyle\scriptstyle{#1}}%
{\XXint\scriptstyle\scriptscriptstyle{#1}}%
{\XXint\scriptscriptstyle\scriptscriptstyle{#1}}%
\!\int}
\def\XXint#1#2#3{{\setbox0=\hbox{$#1{#2#3}{\int}$ }
\vcenter{\hbox{$#2#3$ }}\kern-.6\wd0}}

\def\dashint{\Xint-}

\DeclareMathOperator{\sinc}{sinc}

\def\T2{{\mathbb T}^2}

\def\T2{{\mathbb T}^2}
\def\T3{{\mathbb T}^3}

\newcommand{\ma}{>}
\newcommand{\ap}{``}
\newcommand{\mi}{<}

\newcommand{\s}{\hspace{4ex}}

\newcommand{\ac}{\`}
\newcommand{\ep}{\varepsilon}

\newcommand{\finedim}{{\unskip\nobreak\hfil\penalty50
   \hskip2em\hbox{}\nobreak\hfil\mbox{\rule{1ex}{1ex} \qquad}
   \parfillskip=0pt \finalhyphendemerits=0\par\medskip}}

\begin{document}
\title{
{\bf\Large Approximation results for a general class of  Kantorovich type operators}}
\author{{\bf\large Gianluca Vinti}
\hspace{2mm}
{\bf\large Luca Zampogni}\vspace{1mm}\\
{\it\small Dipartimento di Matematica e Informatica}\\ {\it\small Universit\ac a degli Studi di Perugia},
{\it\small Via Vanvitelli 1, 06123, Perugia (Italy)}\\
{\it\small e-mail:gianluca.vinti@unipg.it, luca.zampogni@unipg.it}}\vspace{1mm}
\maketitle
\begin{center}
{\bf\small Abstract}
\vspace{3mm}
\hspace{.05in}\parbox{4.5in}
{{\small We introduce  and study a family of integral  operators in the  Kantorovich sense  for functions acting on locally compact topological groups.  
We obtain convergence results for the above operators with respect to the pointwise and uniform convergence and in the setting of Orlicz spaces with respect to the modular convergence. Moreover, we show  how our theory applies to several classes of integral and discrete operators, as the  sampling, convolution and  Mellin  type operators in the Kantorovich sense, thus obtaining a simultaneous approach for discrete and integral operators. Further, we derive our general convergence results for particular cases of Orlicz spaces, as $L^p-$spaces, interpolation spaces and exponential spaces. Finally we construct some concrete example of our operators and we show some graphical representations.
}}
\end{center}

 \section{Introduction}
\def\theequation{1.\arabic{equation}}\makeatother
\setcounter{equation}{0}
The goal of the present  paper  is to introduce  and study a family of   operators \ap \`a la Kantorovich'' for functions acting on locally compact topological groups. These operators include, in particular, the family of Kantorovich sampling type operators, which play a crucial role in the theory of Signal and Image Processing.   

Kantorovich sampling type operators were introduced in \cite{BBSV3} and studied in \cite{BM1,BM2,VZ1,VZ2,CV1,CV2}; see also \cite{KA,BBSV1}. In \cite{BBSV3} the authors considered operators acting on functions defined on the  real line with the idea of introducing a Kantorovich version of the generalized sampling operators. Here, we extend the theory to the general case when the underlying space is a locally compact topological group. 

There are various reasons for which it is worth considering Kantorovich sampling operators in the setting of topological groups.  First of all, one has the possibility to retrieve, from our operators,  several families of Kantorovich type operators;  second, one has the possibility to treat simultaneously both discrete and integral operators and both one and multidimensional operators. In particular, we will be able to define and study Kantorovich versions of sampling,  convolution and Mellin operators all in one and multidimensional setting.

The starting point to understand the significance of sampling operators of Kantorovich type is the following (see \cite{BBSV3}): for a locally integrable function $f$ defined on $\mathbb R$, we define  
\begin{equation}\label{1i} (S_wf)(x)=\sum_{-\infty}^\infty\chi(wx-k)w\int_{\frac{k}{w}}^{\frac{k+1}{w}}f(s)ds,\;x\in\mathbb R\s(w>0),\end{equation} where $\chi$ is a suitable kernel  and $f$ is chosen in such a way that the series \eqref{1i} converges.

One of the main differences between the series \eqref{1i} and the generalized  sampling series lies in the use of the mean values $\displaystyle{w\int_{\frac{k}{W}}^{\frac{k+1}{w}}f(s)ds}$ instead of the sampling  values $f(k/w)$. This is because, in practical  situations, more information is usually known around a point than exactly at that point, and therefore  Kantorovich sampling series \eqref{1i} arises as a natural modification of the generalized sampling series in order to reduce time jitter errors; concerning sampling type series, the reader can see  the book \cite{BMV} (Chapt. 8 and 9) and \cite{STEN1,BU,RS,BRS,BSPS,BH,BS1,BS2,BFS,BV1,BV2,BL1,STEN2,BL2,BSS,VI1,BV3,AV1,AV2,VI2,BBSV2,
BBSV4}). Moreover, by using the series \eqref{1i}, one can deal with integrable  (e.g. discontinuous) functions as well, and in fact with functions lying in a general  Orlicz space. This is another important difference between \eqref{1i} and the generalized sampling series, since the  infinite sum of the latter is  not suitable  for integrable functions because it depends on function values  $f(k/w).$  In addition, the generalized sampling operators are not continuous in  $L^p(\mathbb R),$ while operators  \eqref{1i} are instead continuous in $L^p(\mathbb R)$. 

In this paper, we define operators which are analogous to \eqref{1i}, but where both the function $f$ and the sample values $k/w$ are defined in locally compact   topological groups with regular Haar measures. In order to deal with   this general setting,  we have to introduce a whole framework in which a series like \eqref{1i} makes sense and can be studied in detail. 
Namely,  we introduce the family $(S_w)_{w>0}$  of  integral operators defined as follows:
\begin{equation}\label{2i} S_wf(z)=\int_H\chi_w(z-h_w(t))\left[\dfrac{1}{\mu_G(B_w(t))}\int_{B_w(t)}f(u)d\mu_G(u)\right]d\mu_H(t),
\end{equation}
where $f:G \rightarrow \mathbb R$ is a measurable function for which the above integrals are well defined. 
\\
The above operators \eqref{2i} contains, as particular cases, the Kantorovich sampling, convolution and Mellin  type operators, as we will show along the paper. 

One of the main results obtained reads as follows: 
\begin{equation*}
\lim_{w\rightarrow\infty}||S_wf-f||_\infty=0                                                                                                                                                                        \end{equation*}
 for  $f\in C(G)$, where  $C(G)$ is the space of uniformly continuous and bounded functions. 

In the case of not necessarily continuous functions, we obtain a modular convergence result in a suitable subset $\mathcal Y$ of an Orlicz space, i.e. we prove that there exists a constant $\lambda\ma0$ such that $$\lim_{w\rightarrow\infty}I^G_\varphi[\lambda(S_wf-f)]=0,$$ 
for $f\in\mathcal Y$ and where $\varphi$ is a convex $\varphi$-function. 
\\
As pointed out, in the particular cases of operators of Kantorovich type, $\mathcal Y$  coincides with the whole Orlicz space $L^{\varphi}.$

Moreover, we show how the general setting of Orlicz spaces allows us to apply our theory to well known examples of function spaces, as the $L^p-$spaces, the interpolations spaces (or Zygmund spaces) and the exponential spaces, the last two being very useful for PDE's and for embedding theorems. Moreover,  the last example concerning  exponential spaces is of  particular interest since, in this context, the norm convergence is not equivalent to the modular one. 

Finally we construct concrete examples of operators, by using kernels both with compact       and without compact support, and we show with some graphical examples  how the considered  operators reconstruct the function in various situations.

\section{Preliminaries and notations}
\def\theequation{2.\arabic{equation}}\makeatother
\setcounter{equation}{0}
\noindent This section provides the background material which is needed throughout the paper.

First we review some notions concerning topological groups.  In the following, we will  deal with  locally compact Hausdorff topological groups equipped with regular measures. If $H$ is such a topological group, we will denote by $\theta_H$ its neutral element. For simplicity, we will denote the group operation in $H$ by the symbol $+$. It is well-known (see, e.g. \cite{M,PON,BREN}) that there exists a unique (up to multiplication by a constant) left (resp. right) translation-invariant regular Borel measure $\mu_H$ (resp. $\nu_H$). In general, $\mu_H$ and $\nu_H$ are different, and indeed, if $A\subset H$ is a Borel set and if $-A$ denotes the set $-A:=\{-a\;|\;a\in A\}$, there exists a constant $\kappa\ma0$ such that $\mu(-A)=\kappa \nu(A).$ The measures $\mu_H$ and $\nu_H$ coincide  if and only if $H$ is an \emph{unimodular} group.
Examples of unimodular groups are abelian groups, compact groups and discrete groups. In unimodular groups hence the right and the left invariant measures coincide, and we will denote simply by $\mu_H$ this measure. Note that, for unimodular groups, one has $\mu_H(A)=\mu_H(-A)$ for every Borel set $A\subset H$.

Let us now consider a locally compact abelian topological group $G$ with neutral element $\theta_G$. Then  $G$ is an unimodular group.  It can be shown (see \cite{PON}) that a countable symmetric local base $\mathcal B$ of the neutral element $\theta_G$ can be chosen in such a way that
\begin{itemize}

\item[($	\ast$)] if $U\in\mathcal B$ then there exists $V\in\mathcal B$ such that $V+V\subset U$ (and $V-V\subset U$).
\end{itemize}
From now on, when we will make use of a local base $\mathcal B$ of the neutral element $\theta_G$ of an abelian topological group $G$, we will agree that $\mathcal B$ satisfies the condition ($\ast$).

We now move our attention to some basic results on Orlicz spaces.
Let $\varphi:\mathbb R_0^+\rightarrow\mathbb R_0^+$ be a continuous function. We say that $\varphi$ is a \emph{$\varphi$-function} if moreover:
\begin{itemize}
 \item[(a)] $\varphi(0)=0$ and $\varphi(u)\ma0$ for all $u\ma0$;
\item[(b)] $\varphi$ is non-decreasing on $\mathbb R_0^+$;
\item[(c)] $\displaystyle{\lim_{u\rightarrow\infty}\varphi(u)=+\infty}$.
\end{itemize}

Let $G$ be a locally compact abelian topological group with regular Haar measure $\mu_G$.
Let us denote by $M(G)$ the set of measurable bounded functions $f:G\rightarrow\mathbb R$.
Further, by $C(G)$ (resp. $C_c(G)$) we denote the set of functions $f:G\rightarrow\mathbb R$ which are
uniformly continuous and bounded (resp. continuous and with compact support), equipped with the standard $||\cdot||_\infty$ norm, where uniform continuity on $G$ means that:
 for every $\ep\ma0$, there exists a compact set $B_\ep\in\mathcal B$ such that for every $s,z \in G$ with $s-z \in B_\ep,$ then $|f(s)-f(z)|\mi\ep.$   

If a $\varphi$-function $\varphi$ is chosen, one can define a functional $I_\varphi^G:M(G)\rightarrow[0,\infty]$ by $$I^G_\varphi(f):=\int_G\varphi(|f(t)|)d\mu_G(t).$$ $I^G_\varphi$ is a \emph{modular functional} on $M(G)$: it generates the Orlicz space $$L^\varphi(G):=\{f\in M(G)\;|\;I^G_\varphi(\lambda f)\mi\infty\;\mbox{for some}\;\lambda\ma0\}.$$ The subset of $L^\varphi(G)$ consisting of those $f\in M(G)$ for which $I^G_\varphi(\lambda f)\mi\infty$ \emph{for every} $\lambda\ma0$ is denoted by $E^\varphi(G)$. In general one has $E^\varphi(G)\subset L^\varphi(G)$, and they coincide if and only if $\varphi$ satisfies the so called \emph{$\Delta_2$-condition}, i.e.,
\begin{center}\it  there exists a number $M\ma0$ such that $\dfrac{\varphi(2u)}{\varphi(u)}\leq M$ for every $u\ma0.$.\end{center}

There are two different kinds of convergence which are usually used in the context of Orlicz spaces. The first one is determined by a norm on $L^\varphi(G),$ called the \emph{Luxemburg norm} and defined as $$||f||_\varphi:=\inf\{\lambda\ma0\;|\;I^G_\varphi(f/\lambda)\leq \lambda\}.$$ The second is a weaker kind of convergence, called \emph{modular convergence}:

\it a sequence $(f_n)\subset L^\varphi(G)$ converges modularly to $f\in L^\varphi(G)$ if  $$\lim_{n\rightarrow\infty}I^G_\varphi[\lambda(f_n-f)]=0,$$ for some $\lambda\ma0$. \rm

Clearly if $(f_n)_n\subset L^\varphi(G)$ converges in the Luxemburg norm to $f$, then it converges modularly to $f$ as well. The converse is true if and only if the $\Delta_2$-condition is satisfied by $\varphi$.

Orlicz spaces  are natural generalizations of $L^p$ spaces, and in fact if $1\leq p\mi\infty$ and $\varphi(u)=u^p$, the Orlicz space generated by $I^G_\varphi$ is exactly the Lebesgue space $L^p(G)$. The function $\varphi(u)=u^p$ satisfies the $\Delta_2$-condition, hence Luxemburg and modular convergences are the same in $L^p(G)$ and coincide with the convergence with respect to  the standard  norm.

There are other examples of Orlicz spaces which play an important role  in  functional analysis and PDEs. For instance, if we set $\varphi_{\alpha,\beta}(u)=u^\alpha\ln^\beta(e+u)$ ($\alpha\geq 1,\;\beta\ma0$), we obtain the \emph{interpolation space} (also called \emph{$L^\alpha\log^\beta L$-space}) $L^{\varphi_{\alpha,\beta}}(G);$ (see e.g \cite{STEI1,STEI2}.

As another example,
if $\alpha\ma0$, we can take $\varphi_\alpha(u)=\exp\left(u^\alpha\right)-1 $ ($u\in\mathbb R_0^+$). The Orlicz space obtained via $\varphi_\alpha$ is called the \emph{exponential space} $L^{\varphi_\alpha}(G)$ (see \cite{HEN}). This last example is particularly interesting because the function $\varphi_{\alpha,\beta}$ does not satisfy the $\Delta_2$-condition, hence in the space $L^{\varphi_{\alpha,\beta}}(G)$ Luxemburg and modular convergence are different.

For more information on Orlicz spaces and related topics, the reader can be addressed to \cite{KR,MU,RR1,RR2,BMV}.

\section{Approximation results}
\def\theequation{3.\arabic{equation}}\makeatother
\setcounter{equation}{0}
Let $H$ and $G$ be  locally compact Hausdorff topological groups with regular Haar measures $\mu_H$ and $\mu_G$ respectively.  Let us denote by $\theta_H$ (resp. $\theta_G$) the neutral element of $H$ (resp. $G$). We further assume that $G$ is abelian. Let  $\mathcal B\subset G$ be a countable local base of the neutral element $\theta_G$ (which satisfies condition ($\ast$) of the previous section), ordered  by inclusion.
For every $w\ma0$, let $h_w:H\rightarrow G$ be a  map which restricts to a homeomorphism from $H$ to $h_w(H)$.

Let us further assume that for every $w\ma0$, there exists a family $\mathcal B_w=(B_w(t))_{t\in H}$ $\subset G$ of open nonempty subsets of $G$
such that
\begin{itemize}
\item[(i)] $0\mi\mu_G(B_w(t))\mi\infty$ for every $t\in H$ and $w\ma0$;
\item[(ii)] for every $w\ma0$ and $t\in H$, $h_w(t)\in B_w(t)$.
\item[(iii)] if $B\in\mathcal B$, there exists a number $\overline w\ma0$ such that for every $w\ma\overline w$ we have $h_w(t)-B_w(t)\subset B$, for every $t\in H$.
\end{itemize}

Let $(\chi_w)_{w\ma0}$ be a  family of measurable kernel functionals; i.e., $\chi_w:G\rightarrow\mathbb R$,     $\chi_w\in L^1(G)$ and is bounded in a neighborhood of $\theta_G$ ($w\ma0$). We assume that

\begin{itemize}
 \item[($\chi_1$)] the map $t\mapsto \chi_w(z-h_w(t))\in L^1(H)$ for every $z\in G$;
\item[($\chi_2$)] for every $w\ma0$ and $z\in G$, $$\int_H\chi_w(z-h_w(t))d\mu_H(t)=1;$$
\item[($\chi_3$)] for every $w\ma 0$, $$m_{0.\pi}(\chi_w):=\sup_{z\in G}\int_H\left|\chi_w(z-h_w(t))\right|d\mu_H(t)\mi M \mi+\infty;$$
\item[($\chi_4$)] if $w\ma0$, $z\in G$  and $B\in\mathcal B$, set $B_{z,w}=\{t\in H\;|\;z-h_w(t)\in B\}\subset H$. Then
$$\lim_{w\rightarrow\infty}\int_{H\setminus B_{z,w}}\left|\chi_w(z-h_w(t))\right|d\mu_H(t)=0$$ uniformly with respect to $z\in G$;
\item[($\chi_5$)] for every $\ep\ma0$ and compact set $K\subset G$, there exists a symmetric compact set $C\subset G$ containing $\theta_G$ with $\mu_G(C)\mi\infty$ and such that $$\int_{z\notin C}\Upsilon_w(K)\left|\chi_w(z-h_w(t))\right|d\mu_G(z)\mi\ep$$ for every sufficiently large $w\ma 0$ and $h_w(t)\in K$, where $$\Upsilon_w(K):=\mu_H\{t\in H\;|\;h_w(t)\in K\}\s(w\ma0).$$
\end{itemize}
If $w\ma0$, we study  the family of operators $\{S_w:M(G)\rightarrow \mathbb R\}_w$ defined as 
\begin{equation}\label{1} S_wf(z)=\int_H\chi_w(z-h_w(t))\left[\dfrac{1}{\mu_G(B_w(t))}\int_{B_w(t)}f(u)d\mu_G(u)\right]d\mu_H(t),
\end{equation}
where $f:G \rightarrow \mathbb R$ is a measurable function such that the above integrals are well defined. 

We make some concrete examples of operators of the kind \eqref{1}.  In Section 4, we will study in detail these examples.
\\
{\bf Kantorovich Sampling Type Operators}. If $H=\mathbb Z$  and  $G=\mathbb R$, we can choose $h_w:\mathbb Z\rightarrow\mathbb R:k\mapsto t_k/w$, where $\{t_k\}_k$ is a sequence of real numbers such that: $(i)$ $t_k\mi t_{k+1}$ ($k\in\mathbb Z$); $(ii)$ there exist numbers $0\mi\delta\mi\Delta$ such that $\delta\mi t_{k+1}-t_k\mi \Delta$ and $B_w(k)=[t_k/w,t_{k+1}/w]$ . We obtain $$S^{(1)}_wf(z)=\sum_{k\in\mathbb Z}\chi_w(z-t_k/w)\left(\dfrac{w}{t_{k+1}-t_k}\int_{t_k/w}^{t_{k+1}/w}f(u)du\right).$$ See e.g. \cite{BBSV3,VZ1,VZ2,CV1,CV2} for a detailed study of these operators.
\\
{\bf Kantorovich Convolution Type Operators}. If $H=G=\mathbb R$ and $h_w(t)=t/w$, we may choose $B_w(t)=[(t-1)/w,(t+1)/w]$ ($w\ma0$).  Then we obtain  $$S_w^{(2)}f(z)=\int_{-\infty}^{\infty}
\chi_w(z-t/w)\left(\dfrac{w}{2}\int_{(t-1)/w}^{(t+1)/w}f(u)du\right)dt.$$ 
\\
Moreover, if we choose $h_w(t)=t$ and $B_w(t)=[t-1/w,t+1/w]$ ($w\ma0$), then we have
$$S_w^{(3)}f(z)=\int_{-\infty}^\infty\chi_w(z-t)\left(\dfrac{w}{2}\int_{t-1/w}^{t+1/w}f(u)du\right)dt.$$
For the theory of classical convolution operators, see e.g.  \cite{BN}.
\\
{\bf Kantorovich Mellin Type Operators}. If $H=G=\mathbb R^+$, then $\mu_H=\mu_G$ is the logarithmic measure, and the group operation is the product.
If $h_w(t)=t$, we take  $B_w(t)=\left[t\dfrac{w}{w+1},t\dfrac{w+1}{w}\right]$ ($w\ma0$) and we have $$S_w^{(4)}f(z)=\int_0^\infty\chi_w(z/t)\left(\dfrac{1}{2\ln(1+1/w)}\int_{t\frac{w}{w+1}}^{t\frac{w+1}{w}}f(u)\dfrac{du}{u}\right)
\dfrac{dt}{t}.$$ We address the reader to the book \cite{M} and \cite{BJ1,BJ2,BM0,BM3} for the theory and results concerning Mellin type operators.\finedim

\bigskip

We move our attention to some preliminary results concerning the structure and the properties of the operators \eqref{1}.
\\
First of all, we observe that the operators $S_w$ map $L^\infty(G)$ into $L^\infty(G)$. In fact,  for every $f \in L^\infty(G),$
\begin{equation*}\begin{split}
|S_wf(z)|\leq &\int_H\Big|\chi_w(z-h_w(t))\Big|\left[\dfrac{1}{\mu_G(B_w(t))}\int_{B_w(t)}|f(u)|d\mu_G(u)\right]d\mu_H(t)\leq\\&\leq ||f||_{\infty}\int_H|\chi_w(z-h_w(t))d\mu_H(t)\leq ||f||_\infty m_{0,\pi}(\chi_w)\end{split}\end{equation*} for $z\in G$, hence $||S_wf||_\infty\leq ||f||_\infty m_{0,\pi}(\chi_w)$ for every $w\ma0$.

From now on, if $A\subseteq G$ is a measurable set with $\mu_G(A)\mi\infty$ and $f:A\rightarrow\mathbb R$  is an integrable function, we will write $$\dashint_A f(u)d\mu_G(u):=\dfrac{1}{\mu_G(A)}\int_Af(u)d\mu_G(u).$$
Next we prove a first result of convergence.

\begin{theorem}\label{p1} Let $f\in C(G)$. Then
\begin{equation*}
\lim_{w\rightarrow\infty}||S_wf-f||_\infty=0.                                                                                                                                                                     \end{equation*}

\end{theorem}
Clearly, Theorem \ref{p1} implies that if $f$ is continuous at a point $z\in G$, then $S_wf(z)$ converges to $f(z)$.

\noindent 
\proof  By the uniform continuity of $f$, for every $\ep\ma0$, there exists a compact set $B_\ep\in\mathcal B$ such that $|f(z)-f(u)|\mi\ep$ whenever $z-u\in B_\ep$. We can choose an open set $B^{(1)}\subset\mathcal B$ such that  $B^{(1)}+B^{(1)}\subset B_\ep$. By (iii), there exists $\overline w\ma0$ such that if $w\ma\overline w$ then $h_w(t)-B_w(t)\subset B^{(1)}$, for every $t\in H$. We can write
\begin{equation*}\begin{split}|S_wf(z)-&f(z)|=\left|
\int_H\chi_w(z-h_w(t))\left[
\dashint_{B_w(t)}(f(u)-f(z))d\mu_G(u)\right]
d\mu_H(t)\right|\leq\\&\leq\int_H|\chi_w(z-h_w(t))|\left[
\dashint_{B_w(t)}|f(u)-f(z)|d\mu_G(u)\right]d\mu_H(t)\leq\\&\leq\int_{B_{z,w}^{(1)}}|\chi_w(z-h_w(t))|\left[
\dashint_{B_w(t)}|f(u)-f(z)|d\mu_G(u)\right]d\mu_H(t)+\\&+
\int_{H\setminus B_{z,w}^{(1)}}|\chi_w(z-h_w(t))|\left[
\dashint_{B_w(t)}|f(u)-f(z)|d\mu_G(u)\right]d\mu_H(t)=\\&=:I_1+I_2.
\end{split}\end{equation*}
We estimate $I_1$. If $z-h_w(t)\in B^{(1)}$, then, since $h_w(t)-B_w(t)\in B^{(1)}$  for $w\ma\overline w,$ we have $z-u=z-h_w(t)+h_w(t)-u\in B^{(1)}+B^{(1)}\subset B_\ep$ whenever $u\in B_w(t),$ for $w\ma\overline w.$ It follows that, for  $w\ma\overline w,$
$$I_1\leq \ep \int_{B_{z,w}^{(1)}}|\chi_w(z-h_w(t))|d\mu_H(t)\leq \ep m_{0,\pi}(\chi_w)\mi \ep M.$$

For $I_2$, we have
$$I_2\leq ||f||_\infty\int_{H\setminus B_{z,w}^{(1)}}|\chi_w(z-h_w(t))|d\mu_H(t).$$ By the assumption ($\chi_4$), $I_2\rightarrow 0$ as $w\rightarrow\infty$, uniformly with respect to $z\in G$. Combining the estimates for $I_1$ and $I_2$ the proof follows at once.\finedim
\noindent

We now turn to a first result of convergence in Orlicz spaces.

\begin{theorem} \label{t2}Let $\varphi:\mathbb R_0^+\rightarrow\mathbb R_0^+$ be a convex $\varphi$-function and let $f\in C_c(G)$. Then $$\lim_{w\rightarrow\infty}||S_wf-f||_\varphi=0.$$\end{theorem}
\noindent
\proof We must show that
 $$\lim_{w\rightarrow\infty}\int_G\left[\varphi(\lambda|S_wf(z)-f(z)|)\right]d\mu_G(z)=0,$$ for every $\lambda\ma0$.

 Let us consider the family of functions $g_w(z):=\varphi(\lambda|S_wf(z)-f(z)|)$. Then $g_w:G\rightarrow\mathbb R$ is nonnegative for every $w\ma0$. Moreover, $\displaystyle{\lim_{w\rightarrow\infty}g_w(z)=0}$ uniformly with respect to $z\in G$. This is due to the  the properties of $\varphi$ and the fact that $||S_wf-f||_\infty\rightarrow 0$ as $w\rightarrow\infty$ (see Theorem \ref{p1}). Now we show that it is possible to apply the Vitali convergence Theorem to the family $(g_w)_w$.

Let $K_1=Supp(f).$ Choose a symmetric compact set $K\subset G$  satisfying $K_1\varsubsetneq K$.  By (iii),  for sufficiently large $w\ma0$, $B_w(t)\subseteq h_w(t)+B$ ($B\in\mathcal B$). Then $B_w(t)\cap K_1=\emptyset$ for sufficiently large $w\ma0$ and therefore if $h_w(t)\notin K$ we have $$\int_{B_w(t)}f(u)du=0.$$

Now, fix  $\lambda\ma0$, let $\ep\ma0$ and let $C\subset G$ with $\mu_G(C)\mi\infty$ be such that ($\chi_5$) is valid for $K$. Let us estimate $$I:=\int_{G\setminus C}g_w(z)d\mu_G(z).$$ We use the notation $K_{t,w}=\{t\in H\;|\;h_w(t)\in K\}$ ($t\in H,\;w\ma0$). We have
\begin{equation*}\begin{split}I&=\int_{G\setminus C}\varphi\left[\lambda|S_wf(z)-f(z)|\right]d\mu_G(z)\leq\\&\leq\int_{G\setminus C}\varphi\left(2\lambda||f||_\infty\int_{K_{t,w}}\left|\chi_w(z-h_w(t))\right|d\mu_H(t)\right)d\mu_G(z)\leq\\&\int_{G\setminus C}\left(\dfrac{1}{\Upsilon_w(K)M}\int_{K_{t,w}}\varphi\left(2\lambda M||f||_\infty\right)\Upsilon_w(K)|\chi_w(z-h_w(t))|d\mu_H(t)\right)d\mu_G(z)\leq\\&
\dfrac{1}{\Upsilon_w(K)M}\int_{K_{t,w}}\left(\varphi(2\lambda M||f||_\infty )\int_{G\setminus C}\Upsilon_w(K)|\chi_w(z-h_w(t))|d\mu_G(z)\right)d\mu_H(t)\leq\\&\dfrac{\ep}{\Upsilon_w(K)M}\int_{K_{t,w}}\varphi(2\lambda M||f||_\infty)d\mu_H(t)=\dfrac{\ep\varphi(2\lambda M||f||_\infty)}{M}\mi\infty,
\end{split}\end{equation*}
where we used the Jensen's inequality, the Fubini-Tonelli Theorem and ($\chi_5$).
Moreover, it is easy to see that, for every measurable set $A\subset G$ with $\mu_G(A)\mi\infty$, we have
$$\int_A\varphi(\lambda|S_wf(z)-f(z)|)d\mu_G(z)\leq\int_A\varphi(2\lambda M||f||_\infty)d\mu_H(t)=\varphi(2\lambda M||f||_\infty)\mu_G(A).$$ So, for fixed $\ep\ma0$, it suffices to take $\delta\mi\dfrac{\ep}{\varphi(2\lambda M||f||_\infty)}$, to obtain $$\int_A \varphi(\lambda|S_wf(z)-f(z)|)d\mu_G(z)\leq\ep,$$ for every measurable set $A\subset G$ with $\mu_G(A)\mi\delta.$ The Vitali convergence theorem can be applied, and the theorem is therefore proved. \finedim
\noindent

It is a matter of fact that, except for the standard operators defined on $\mathbb R$, when one has to face the problem of the convergence in the space $L^\varphi(G)$ for a generic function $f\in L^\varphi(G)$, then one has to make one additional assumption, which allows to compare the value of an integral over $H$ of the function $\varphi$ which involves $f$ and the sets $B_w(t)$  (in a sense specified below) with the value $I^G_\varphi(\lambda f)$. We formulate this assumption:
\begin{itemize}
\item[($\chi_6$)] We assume that there exists a vector subspace $\mathcal Y\subset L^\varphi(G)$ with $C_c^\infty(G)\subset \mathcal Y$  and such that, for every $g\in\mathcal Y$, there holds

\begin{equation}\label{cond}\limsup_{w\rightarrow\infty}||\chi_w||_{L^1(G)}  I^H_\varphi \left(\dashint_{B_w(\cdot)}g(z)d\mu_G(z)\right)\leq CI^G_\varphi( g),\end{equation} for some $C\ma0$.\end{itemize}
We remark that condition \eqref{cond} is not assured in general: however,
 in the cases when $H$ and $G$ are subgroups of $\mathbb R$, as we will see below,  the \ap Kantorovich'' nature of  the operators considered here allows us to discharge ($\chi_6$) on the kernels $\chi_w$.

Under this additional assumption we can prove the following
\begin{theorem}\label{t3}Let $\varphi$ be a convex $\varphi$-function and let $f\in\mathcal Y$, where $\mathcal Y$ is defined in ($\chi_6$). Then, if $\lambda\ma0$,
\begin{equation*}
I^G_\varphi(\lambda S_wf)\leq\dfrac{C}{M}I^\varphi_G(\lambda Mf)
\end{equation*}
for sufficiently large $w\ma0$.

In particular $S_w:\mathcal Y\rightarrow L^\varphi(G)$ is well defined for every $w\ma0$.\end{theorem}
\noindent
\proof Let $\lambda\ma 0$ be such that the quantity $I^G_\varphi[\lambda Mf]\mi\infty$. Then, using \eqref{cond} with $g=\lambda Mf$,  \begin{equation*}\begin{split}
&I^G_\varphi[\lambda S_wf]=\int_G\varphi(\lambda |S_wf(z)|)d\mu_G(z)\leq \\&\leq
\int_G\varphi\left[\lambda\int_H|\chi_w(z-h_w(t))|
\left(\left|\dashint_{B_w(t)}f(u)d\mu_G(u)\right|\right)d\mu_H(t)\right]d\mu_G(z)
\leq\\&\leq\dfrac{1}{M}\int_G\left[\int_H\varphi\left(\lambda M\left|\dashint_{B_w(t)}f(u)
d\mu_G(u)\right|\right)|\chi_w(z-h_w(t))|d\mu_H(t)\right]d\mu_G(z)
\leq\\&\leq\dfrac{1}{M}\left[\int_H\varphi\left(\lambda M\left|\dashint_{B_w(t)}f(u)d\mu_G(u)\right|\right)\left(\int_G|\chi_w(z-h_w(t))|d\mu_G(z)\right)
\right]d\mu_H(t)\leq\\&\leq \dfrac{1}{M}||\chi_w||_{L^1(G)} I^H_\varphi\left(\lambda M\dashint_{B_w(\cdot)}f(u)d\mu_G(u)\right)\leq \dfrac{C}{M}I^G_\varphi[\lambda Mf],
\end{split}\end{equation*}
for sufficiently large $w\ma0$. \finedim
\noindent

Our next result concerns the convergence in the Orlicz space $L^\varphi(G)$ of $S_wf$ to $f$ as $w\rightarrow\infty$. To prove it, we need the following lemma (see \cite{BM00,BV1,BMV}).
\begin{lemma}The set $C_c^\infty(G)$ is dense in $L^\varphi(G)$ with respect to the modular convergence.\end{lemma}

Now, the main convergence result can be stated as follows:
\begin{theorem}\label{t6} Let $\varphi$ be a convex $\varphi$-function and let $f\in\mathcal Y$. Then there exists a constant $\lambda\ma0$ such that $$\lim_{w\rightarrow\infty}I^G_\varphi[\lambda(S_wf-f)]=0.$$\end{theorem}
\noindent
\proof Let $\ep\ma0$. By the above Lemma, we can find a function $g\in C_c(G)$ and a constant $\eta\ma0$ such that
$$I^G_\varphi[\eta(f-g)]\mi\ep.$$
Choose $\lambda\ma0$ such that $3\lambda(1+M)\mi\eta.$
We can write
\begin{equation*}\begin{split}I^G_\varphi&[\lambda(S_wf-f)]\leq
I^G_\varphi[3\lambda(S_wf-S_wg)]+I^G_\varphi[3\lambda(S_wg-g)]+
I^G_\varphi[3\lambda(f-g)]\leq\\&\leq \dfrac{C}{M}I^G_\varphi[\eta (f-g)]+I^G_\varphi[3\lambda(S_wg-g)]+I^G_\varphi[\eta(f-g)]\leq \\&\leq \left(\dfrac{C}{M}+1\right)\ep+I^G_\varphi[3\lambda(S_wg-g)].
\end{split}\end{equation*}
From Theorem \ref{t2} the proof follows easily since $\ep$ is arbitrarily chosen. \finedim
\noindent
\section{Applications}
\def\theequation{4.\arabic{equation}}\makeatother
\setcounter{equation}{0}

In this section, we will give concrete examples of applications of the theory developed in the previous section.
Some of these examples are known in the literature, whereas others are generalizations of well known  operators to the \ap Kantorovich'' setting.

\begin{itemize}
 \item[\bf(1)] We begin with a kind of operators discussed in \cite{BBSV3}. Let $H=(\mathbb Z,+)$ and $G=(\mathbb R,+)$ provided with the counting and Lebesgue measures respectively. Let us define $h_w:\mathbb Z\rightarrow\mathbb R:k\mapsto t_k/w$ ($w\ma0$), where $t_k\mi t_{k+1}$ for every $k\in\mathbb Z$, and $\delta\mi t_{k+1}-t_k\mi\Delta$ for some numbers $0\mi\delta\mi\Delta\mi\infty$. Set $\Delta_k=t_{k+1}-t_k.$ ($k\in\mathbb Z$), and let $B_w(k)=[t_k/w,t_{k+1}/w].$ If $f:\mathbb R\rightarrow\mathbb R$ is a measurable function, the corresponding family of operators is defined as
$$S_w^{(1)}f(z)=\sum_{k\in\mathbb Z}\chi_w(z-t_k/w)\dfrac{w}{\Delta_k}\int_{t_k/w}^{t_{k+1}/w}f(u)du.$$

For the family $(S_w^{(1)})_{w\ma0},$ a theory has already been introduced in \cite{BBSV3}. We show that the assumptions ($\chi_1$)--($\chi_6$) reduce to those in the above paper.

Assumptions from ($\chi_1$) to ($\chi_3$) can be easily rewritten with $\displaystyle{\sum_{k\in\mathbb Z}}$ instead of $\displaystyle{\int_H\dots d\mu(t)}$.

($\chi_4$) can be rewritten in a more familiar form as follows:
\ap for every $\gamma\ma0$
$$\lim_{w\rightarrow\infty}\sum_{|z-t_k/w|\ma\gamma}|\chi_w(z-t_k/w)|=0,$$ uniformly with respect to $z\in\mathbb R$'';

($\chi_5$) is equivalent to the following:
\ap for every $\ep\ma0$ and $\gamma\ma0$ there exists a number $M\ma0$ such that
$$\int_{|z|\ma M}w|\chi_w(z-t_k/w)|dz\mi\ep$$ for sufficiently large $w\ma0$ and $t_k/w\in[-\gamma,\gamma]$''. Indeed, in this case the compact set $K$ can be taken as $[-\gamma,\gamma]$ ($\gamma\ma0$), while the symmetric compact set $C$ is given by $[-M,M]$. Now, if we compute the quantity $\Upsilon_w([-\gamma,\gamma])$ we have
\begin{equation*}\begin{split}\dfrac{2\gamma w}{\Delta}-2\leq\Upsilon_w([-\gamma,\gamma])\leq 2\left(\dfrac{\gamma w}{\delta}+2\right)=4+2\dfrac{\gamma}{\delta}w.\end{split}\end{equation*}
It follows that
\begin{equation*}\begin{split}\int_{|z|\ma M}\left(\dfrac{2\gamma w}{\Delta}-2\right)|\chi_w(z-t_k/w)|dz&\leq\int_{|z|\ma M}\Upsilon_w([-\gamma,\gamma])|\chi_w(z-t_k/w)|dz\\ &\leq \int_{|z|\ma M}\left(4
+2\dfrac{\gamma w}{\delta} \right)|\chi_w(z-t_k/w)|dz.\end{split}\end{equation*}
 So, since $X_w\in L^1(\mathbb R)$ for every $w\ma0$,  ($\chi_5$) is equivalent to $$\int_{|z|\ma M}w|\chi_w(z-t_k/w)|dz\mi \ep$$ for sufficiently large $w\ma0$ and $t_k/w\in[-\gamma,\gamma]$;

($\chi_6$) translates to a condition on the kernels: indeed \eqref{cond} becomes
\begin{equation*}\begin{split}&||\chi_w||_1 I^\mathbb Z_\varphi\left[\alpha \left(\dfrac{w}{\Delta_k}\int_{t_k/w}^{t_{k+1}/w}f(z)dz\right)\right]\leq\\&\leq ||\chi_w||_1\dfrac{w}{\delta}\left(\sum_{k\in\mathbb Z}\int_{t_k/w}^{t_{k+1}/w}\varphi[\alpha|f(z)|]dz\right)=\\&=||\chi_w||_1\dfrac{w}{\delta} I_\varphi^\mathbb R[\alpha f] .\end{split}\end{equation*}
So it suffices to have $$\limsup_{w\rightarrow\infty}w||\chi_w||_1\mi\infty.$$ This happens, for example, if $\chi_w(s)$ is given by $\chi_w(s)=\chi(ws)$ where $\chi\in L^1(\mathbb R)$. Actually, in \cite{BBSV3}, the kernels under consideration have exactly this form.
Therefore, in this case, we have $\mathcal Y=L^\varphi(\mathbb R)$.

We now slightly modify the operator $S_w^{(1)}$ to show an interesting consequence of the Kantorovich frame. So, let $H=(\mathbb Z,+)$ and $G=(\mathbb R,+)$ provided with the counting and Lebesgue measures respectively as before, and consider the map $h_w:\mathbb Z\rightarrow\mathbb R: k\mapsto t_k,$ where the sequence $(t_k)_k\subset\mathbb R$ is chosen as above.  Now, set $B_w(k):=[t_k-1/w,t_k+1/w]$ ($w\ma0$).
With these choices, we obtain the operators $$S_w^{(1,1)}f(z)=\sum_{k\in\mathbb Z}\chi_w(z-t_k)\dfrac{w}{2}\int_{t_k-1/w}^{t_{k}+1/w}f(u)du.$$ For the family $(S_w^{(1,1)})_w$, we have the same results as for the family $(S_w^{(1)})_w$ above. In particular, ($\chi_6$) translates to a condition on the family of kernels $(\chi_w)_w$ which is satisfied if $\chi_w(s)=\chi(ws)$ where $\chi\in L^1(\mathbb R)$, and so $\mathcal Y=L^\varphi(\mathbb R)$ (see the discussion above). However, this examples has an interesting feature. For fixed $k\in\mathbb Z$, if $f\in L^1_{loc}(\mathbb R)$, the quantity $$\dfrac{w}{2}\int_{t_k-1/w}^{t_{k}+1/w}f(u)du$$ converges, as $w\rightarrow\infty$, to the value $f(t_k)$, for a.a. $t_k$.  This is the Lebesgue-Besicovich Differentiation Theorem (see \cite{EG}). Roughly speaking, this fact tells us that, for large $w\ma0$, the operators $S_w^{(1,1)}f(\cdot)$ can be asymptotically compared with the classical sampling operators $$S_wf(z)= \sum_{k\in\mathbb Z}\chi_w(z-t_k)f(t_k).$$
To retrieve the exact formula of the classical sampling operators, it suffices to take $t_k = s_k/w$ and $\chi_w(z)=\chi(wz).$ Hence we obtain
$$S_wf(z)=\sum_{k\in\mathbb Z}\chi(wz-s_k)f\left(\dfrac{s_k}{w}\right).$$

 However, it is well known (see \cite{BMV}) that the sampling operators $S_wf(\cdot)$ do not converge to $f(\cdot)$ for an arbitrary function $f\in L^\varphi(\mathbb R)$. In fact, one can show that the convergence is assured in a proper subspace $\mathcal Y$ of $L^\varphi(\mathbb R)$, namely $\mathcal Y=BV^\varphi(\mathbb R)\cap E^\varphi(\mathbb R)$ (where $BV^\varphi(\mathbb R)$ is the subset of $M(\mathbb R)$ consisting of those $f\in M(\mathbb R)$ satisfying  $\varphi(\lambda|f|)\in BV(\mathbb R)$ for every $\lambda\ma0$). This shows the importance of  the regularizing property of the mean value in the Kantorovich type operators.
\item[\bf(2)] We take $H=G=\mathbb R$, $h_w:\mathbb R\rightarrow\mathbb R:t\mapsto t/w$ and $B_w(t)=[(t-1)/w,(t+1)/w]$ ($w\ma0$, $t\in\mathbb R$). In this case we have
$$S_w^{(2)}f(z)=\int_{-\infty}^\infty \chi_w(z-t/w)\left(\dfrac{w}{2}\int_{(t-1)/w}^{(t+1)/w}f(s)ds\right)dt,$$
 which is a Kantorovich version of a convolution operator.
We rewrite the assumptions ($\chi_1$)--($\chi_6$).

Again, assumptions from ($\chi_1$) to ($\chi_3$) may be rewritten with $H=\mathbb R$ and the standard Lebesgue measure on $\mathbb R$ in place of $d\mu(t)$ and $d\mu(z)$.

($\chi_4$) can be easily written as $$\lim_{w\rightarrow\infty}\int_{|z-t/w|\ma\gamma}|\chi_w(z-t/w)|dt=0,$$ for every $\gamma\ma0$ and uniformly with respect to $z\in\mathbb R$;

($\chi_5$) assumes an interesting form: namely, (we can take $K=[-\gamma,\gamma]$ and $C=[-M,M]$)  we require that for every $\ep\ma0$ and $\gamma\ma0$ there exists a number $M\ma0$ such that $$\int_{|z|\ma M} w|\chi_w(z-t/w)|dz\mi\ep,$$ for sufficiently large $w\ma0$ and $t/w\in[-\gamma,\gamma]$. Indeed, in this case he have
$\Upsilon_w([-\gamma.\gamma])=2\gamma w.$

($\chi_6$) can be rewritten by observing that
\begin{equation*}\begin{split}&\s\s||\chi_w||_1\int_\mathbb R\varphi\left[\dfrac{w}{2}\alpha\left|\int_{(t-1)/w}^{(t+1)/w}f(u)du\right|\right]dt\leq \\&\leq ||\chi_w||_1\dfrac{w}{2}\int_0^{2/w}\left(\int_\mathbb R\varphi\left(\alpha|f(s+(t-1)/w)|\right)dt\right)ds,
\end{split}
\end{equation*}
by using the Jensen's inequality, the change of variable $s=u-(t-1)/w$ and the Fubini-Tonelli theorem.
 A further change of variable $\rho=s+(t-1)/w$ gives \begin{equation*}\begin{split}&||\chi_w||_1\int_\mathbb R\varphi\left[\dfrac{w}{2}\alpha\left|\int_{(t-1)/w}^{(t+1)/w}f(u)du\right|\right]dt\leq\\ &\leq ||\chi_w||_1\dfrac{w}{2}\int_0^{w/2}wI^\mathbb R_\varphi(\alpha f)ds=w||\chi_w||_1 I^\mathbb R_\varphi(\alpha f).\end{split}\end{equation*}

Again, the condition ($\chi_6$) is satisfied if $$\limsup_{w\rightarrow\infty}w||\chi_w||_1\mi\infty,$$ which is true if, for instance, $\chi_w(s)=\chi(ws)$ for every $s\in\mathbb R$, where $\chi\in L^1(\mathbb R)$. Again, in this case we have $\mathcal Y=L^\varphi(\mathbb R)$.

\item[\bf(3)] A similar example can be obtained by setting $h_w(t)=t$ for every $w\ma0$ and $t\in\mathbb R.$ In this case, we can take $B_w(t)=[t-1/w,t+1/w],$ and the operators have the form
$$S_w^{(3)}f(z)=\int_\mathbb R\chi_w(z-t)\left(\dfrac{w}{2}\int_{t-1/w}^{t+1/w}f(u)du\right)dt.$$
In this case, however, the statements of assumptions ($\chi_5$) and ($\chi_6$) slightly differ from those mentioned above, namely:

($\chi_5$) $\rightarrow$ for every $\ep\ma0$ and $\gamma\ma0$ there exists a number $M\ma0$ such that $$\int_{|z|\ma M}|\chi_w(z-t)|dz\mi\ep$$ for sufficiently large $w\ma0$ and $t\in[-\gamma,\gamma]$ (indeed, one has $\Upsilon_w([-\gamma,\gamma])=2\gamma$). The above assumption is clearly satisfied if, for example, the family $(\chi_w)_w$ is uniformly bounded by a function $\chi\in L^1(\mathbb R)$.

($\chi_6$) $\rightarrow$ $\displaystyle{\limsup_{w\rightarrow\infty}||\chi_w||_1\mi\infty,}$ arguing as before, and $\mathcal Y=L^\varphi(\mathbb R)$.

Note that, if $f\in L^1_{loc}(\mathbb R)$, in $S_w^{(3)}f(z)$ the factor $\displaystyle{\dfrac{w}{2}\int_{t-1/w}^{t+1/w}f(u)du}$ converges, as $w\rightarrow\infty$, to the value $f(t)$ for a.e. $t\in\mathbb R$ (this is the Lebesgue-Besicovich Differentiation Theorem again). Roughly speaking, this fact tells us that, for large values of $w\ma0$,  $S_w^{(3)}f(z)$ can be related to the standard convolution operator
$$C_wf(z)=\int_\mathbb R\chi_w(z-t)f(t)dt.$$ Compare our results with the very well-known theorem which states that $C_wf(z)\rightarrow f(z)$ as $w\rightarrow\infty,$ when $f\in L^p(\mathbb R^N)$ and $\chi_w$ is a family of mollifiers (any sequence of mollifiers satisfies the  assumptions ($\chi_1$)--($\chi_6$)).

\item[\bf(4)] Let $H=G=\mathbb R^+$. The group operation in $\mathbb R^+$ is the product, hence $\theta_G=1$. Set $h_w(t)=t$ and $B_w(t)=\left[t\dfrac{w}{w+1},t\dfrac{w+1}{w}\right]$ for every $w\ma0$ and $t\ma0$. The only regular Haar measure on $\mathbb R^+$ (up to a multiplicative constant) is the logarithmic measure $\mu(\mathbb R^+)=\displaystyle{\int_{\mathbb R^+} \dfrac{dt}{t}.} $ The family $\mathcal B$ can be taken as $\mathcal B=\{[1/\alpha,\alpha],\;\alpha\ma1\}.$ The Haar measure of $B_w(t)$ is $2\ln\dfrac{w+1}{w}.$ If $f\in M(\mathbb R^+)$, then we obtain
$$S_w^{(4)}f(z)=\int_0^\infty \chi_w\left(\dfrac{z}{t}\right)\dfrac{1}{2\ln(1+1/w)}\left(\int_{t\frac{w}{w+1}}^{t\frac{w+1}{w}} f(u) \dfrac{du}{u}\right) \dfrac{dt}{t}.$$

Assumptions ($\chi_1$)--($\chi_4$) can be easily adapted with $H=G=\mathbb R^+$ and $d\mu(t)=\dfrac{dt}{t}.$

As in the above example, ($\chi_5$) assumes the following form: \ap for every $\ep\ma0$ and $\gamma\ma1$, there exists a number $M\ma1$ such that $$\int_{z\notin[1/M,M]}|\chi_w(z/t)|\dfrac{dz}{z}\mi\ep$$ for sufficiently large $w\ma0$ and $t\in[1/\gamma,\gamma]$'' (indeed, in this case, $\Upsilon([1/\gamma,\gamma])=2\ln\gamma$).

It is possible to prove that  ($\chi_6$) is satisfied if $$\limsup_{w\rightarrow\infty}||\chi_w||_1\mi\infty;$$ one has to keep in mind, however, that the $L^1$-norm is with respect to the measure $\dfrac{dt}{t};$ again, $\mathcal Y=L^\varphi(\mathbb R^+)$.

It is worth to spend again a word concerning the form of $S_w^{(4)}.$ If $w\rightarrow\infty$ and if $f\in L^1_{loc}(\mathbb R^+,dt/t)$, then the value $I:=\displaystyle{\dfrac{1}{2\ln(1+1/w)}\int_{t\frac{w}{w+1}}^{t\frac{w+1}{w}} f(u) \dfrac{du}{u}}$ converges to $f(t)$ for a.e. $t\in\mathbb R^+$. Indeed,  by the change of variables $ts=u$, one gets $$I=\dfrac{1/w}{\ln(1+1/w)}\dfrac{w}{2}\int_{1-\frac{1}{w+1}}^{1+\frac{1}{w}} f(ts) \dfrac{ds}{s}:=\dfrac{1/w}{\ln(1+1/w)}I_1,$$ and $I_1$ converges to $f(t)$ as $w\rightarrow\infty$ for  a.e. $t\in\mathbb R^+$ by  the Lebesgue-Besicovich theorem again, since the term 
$\dfrac{1/w}{\ln(1+1/w)}\dfrac{w}{2}$ is the reciprocal  of the logarithmic (Haar)  measure of the integration set and in fact it is an average.

Now, since the factor $\dfrac{1/w}{\ln(1+1/w)}$ converges to 1 as $w\rightarrow\infty$, then
the operator $S_w^{(4)}f(z)$ can be asymptotically compared, as $w\rightarrow\infty$, to the standard Mellin operator $$M_wf(z)=\int_0^\infty\chi_w\left(\dfrac{z}{t}\right)f(t)\dfrac{dt}{t}.$$
\item[\bf(5)] Our last application concerns operators which approximate functions $f:\mathbb R^N\rightarrow\mathbb R$.  Let $(t_k)_{k\in\mathbb Z}$ be a sequence of real numbers such that:
 $(i)$ $-\infty\mi t_{k}\mi t_{k+1}\mi\infty$; $(ii)$ $\displaystyle{\lim_{k\rightarrow\pm\infty}t_{k}=\pm\infty}$;
$(iii)$ $\delta\mi t_{k+1}-t_{k}\mi\Delta$ for some fixed numbers $0\mi\delta\mi\Delta\mi\infty.$ Let $\Delta_{k}=t_{k+1}-t_{k}.$

 Let $H=\mathbb Z^N$ and $G=\mathbb R^N$. Denote points in $\mathbb Z^N$ by ${\bf k}=(k_1,\dots,k_N).$ For every $w\ma0$, let us define $h_w:\mathbb Z^N\rightarrow\mathbb R^N$ by $h_w({\bf k})=\dfrac{1}{w}{\bf t_{\bf k}}=(t_{k_1}/w,\dots,t_{k_N}/w)\in\mathbb R^N$. Let us consider the grid in $\mathbb R^N$ determined by the point ${\bf t_k}$. Then $R^N$ is divided into parallelepipeds $C_{\bf{k}}=[t_{k_1}/w,t_{k_1+1}/w]\times\dots\times [t_{k_{N}}/w,t_{k_N+1}/w]$. Each parallelepiped $C_{{\bf k}}$ has $N$-dimensional Lebesgue measure equal to $|C_{{\bf k}}|=\dfrac{\Delta_{k_1}\cdots\Delta_{k_N}}{w^N}.$

 If $f\in M(\mathbb R^N)$, we are left with the family
$$S_w^{(5)}f({\bf z})= \sum_{{\bf k}\in\mathbb Z^N}\chi_w({\bf z}-{\bf t_k}/w)\dfrac{1}{|C_{{\bf k}}|}\int_{C_{{\bf k}}}f({\bf u})du.$$

Our theory applies to the family above (see \cite{CV1,CV2,BM1}  for details concerning the series $S_w^{(5)}f$). Moreover, one can define operators analogous to $S^{(2)}_wf,\;S^{(3)}_wf$ and $S^{(4)}_wf$ in a multidimensional setting by adapting the construction we made above to those concrete cases.
\end{itemize}

\begin{remark}\rm  As we mentioned before, in all our examples assumption ($\chi_6$) translates in a condition on the kernels $\chi_w$. This is due to the behaviour of the Kantoro\-vich-type operators in $L^\varphi(\mathbb R^N)$. \end{remark}

\section{Concrete results in Orlicz spaces}
\def\theequation{5.\arabic{equation}}\makeatother
\setcounter{equation}{0}

In this section we rewrite the results obtained before in the general setting to some important cases.

The first case we examine are the standard $L^p$-spaces.  Let  $\varphi(u)=u^p$ ($p\geq 1$), $u\in\mathbb R_0^+$. The Orlicz space $L^\varphi(G)$ coincides with the standard space $L^p(G)$. Theorems \ref{t3} and \ref{t6}  can be summarized as follows:
\begin{theorem}\label{lp} Let $f\in\mathcal Y.$ Then $$||S_wf||_p\leq \left(CM^{p-1}\right)^{1/p}||f||_p.$$ Moreover
$$\lim_{w\rightarrow\infty}||S_wf-f||_p=0.$$

\end{theorem}

Other important examples of Orlicz spaces are furnished by the so-called \ap interpolation spaces'' ($L^\alpha\ln^\beta L$): these spaces are defined as the Orlicz spaces $L^{\alpha,\beta}(G)$, where $\alpha\geq1$ and $\beta\ma 0$. The generating function is given by $\varphi_{\alpha,\beta}(u)=u^\alpha\ln^\beta(e+u)$. Hence $L^{\alpha,\beta}(G)$ is made up of all the functions $f\in M(G)$ such that there exists a constant $\lambda \ma0$ satisfying $$I^{\alpha,\beta}(\lambda f):=\int_\mathbb G\left(\lambda |f(x)|\right)^\alpha \ln^\beta (e+\lambda|f(x)|)dx\mi+\infty.$$
The interested reader can be addressed to \cite{STEI1,STEI2} for more information on interpolation spaces. The function $\varphi_{\alpha,\beta}$ satisfies the $\Delta_2$-condition, hence modular and Luxemburg convergences are the same in $L^{\alpha,\beta}(G)$. Again, we can restate theorems \ref{t3} and \ref{t6} as follows:
\begin{theorem}\label{llog} Let $f\in\mathcal Y.$
Then \begin{equation*}\begin{split}&\int_G|S_wf(x)|^\alpha\ln^\beta\left(e+\lambda|S_wf(x)|\right)dx\leq \\&\leq CM^{\alpha-1}\int_G\left(|f(x)|^\alpha\ln^\beta(e+\lambda M|f(x)|)\right)dx,\end{split}\end{equation*} for every $\lambda\ma0$.  Moreover, $$\lim_{w\rightarrow\infty}||S_wf-f||_{\varphi_{\alpha,\beta}}=0.$$
\end{theorem}

The last example of Orlicz space we consider  is a space where the $\Delta_2$-condition is not fulfilled, hence modular and Luxemburg convergence are distinct. This space is the so called \ap exponential space''. To define an exponential space, let us fix $\alpha\ma0$, and consider the function $\varphi_\alpha:\mathbb R_0^+\rightarrow\mathbb R_0^+:u\mapsto\exp(u^\alpha)-1.$ The space $L^{\varphi_\alpha}(G)$ consists of those functions $f\in M(G)$ for which there exists a constant $\lambda\ma0$ such that $$I^\alpha(\lambda f):=\int_\mathbb G\left(\exp[(\lambda|f(x)|)^\alpha]-1\right)dx\mi +\infty;$$ (for more information on exponential spaces, see, e.g. \cite{HEN}).

Theorems \ref{t3} and \ref{t6} can be stated as
\begin{theorem}\label{exp} Let $f\in \mathcal Y.$
Then for some $\lambda\ma0$ we have

$$\int_G\left(\exp[(\lambda|S_wf(x)|)^\alpha]-1\right)dz\leq \dfrac{C}{M}\int_G
\left(\exp[(\lambda M|f(x)|)^\alpha]-1\right).
$$
Moreover, there exists a number $\lambda\ma0$ such that $$\lim_{w\rightarrow\infty}\int_G\left(\exp[(\lambda|S_wf(x)-f(x)|)^\alpha]-1\right)dx=0.$$
\end{theorem}
\section{Some graphical representations}
\def\theequation{6.\arabic{equation}}\makeatother
\setcounter{equation}{0}

This section provides some graphical representations of the convergence of the operators we have studied in the previous sections. In all the examples below the convergence must be interpreted as to be in the $L^p$ setting.

We will concentrate on the examples (2), (3) and (4) of  Section 4, 
since graphical examples of operators (1) can be found in \cite{BBSV3} and for operators (5) one can see \cite{CV1,CV2}.

Although the prototypical example of kernel is obtained from the Fejer's kernel function $$F(x)=\dfrac{1}{2}\sinc^2\left(\dfrac{x}{2}\right),$$ where $$\sinc(x)=\begin{cases}\dfrac{\sin\pi x}{\pi x},&x\in\mathbb R\setminus\{0\},\\1,&x=0\end{cases},$$
 it will be convenient for computational purposes to take a kernel with compact support over $\mathbb R$. Well known examples of such kernels are those arising from linear combinations the so-called $B$-splines functions of order $n\in\mathbb N$, namely
$$M_n(x)=\dfrac{1}{(n-1)!}\sum_{j=0}^n(-1)^j\begin{pmatrix}n\\j\end{pmatrix}\left(\dfrac{n}{2}+x-1\right)_+^{n-1},$$ where the symbol $(\cdot)_+$ denotes the positive part. Below we represent the graphs of the functions $M_3(x),M_4(x)$ (Figure 1), and $M(x)=4M_3(x)-3M_4(x)$ (Figure 2).
\vspace{-.3cm}
\begin{figure}[htbp]
\begin{center}\includegraphics[width=.55\textwidth,height=.34\textwidth]
{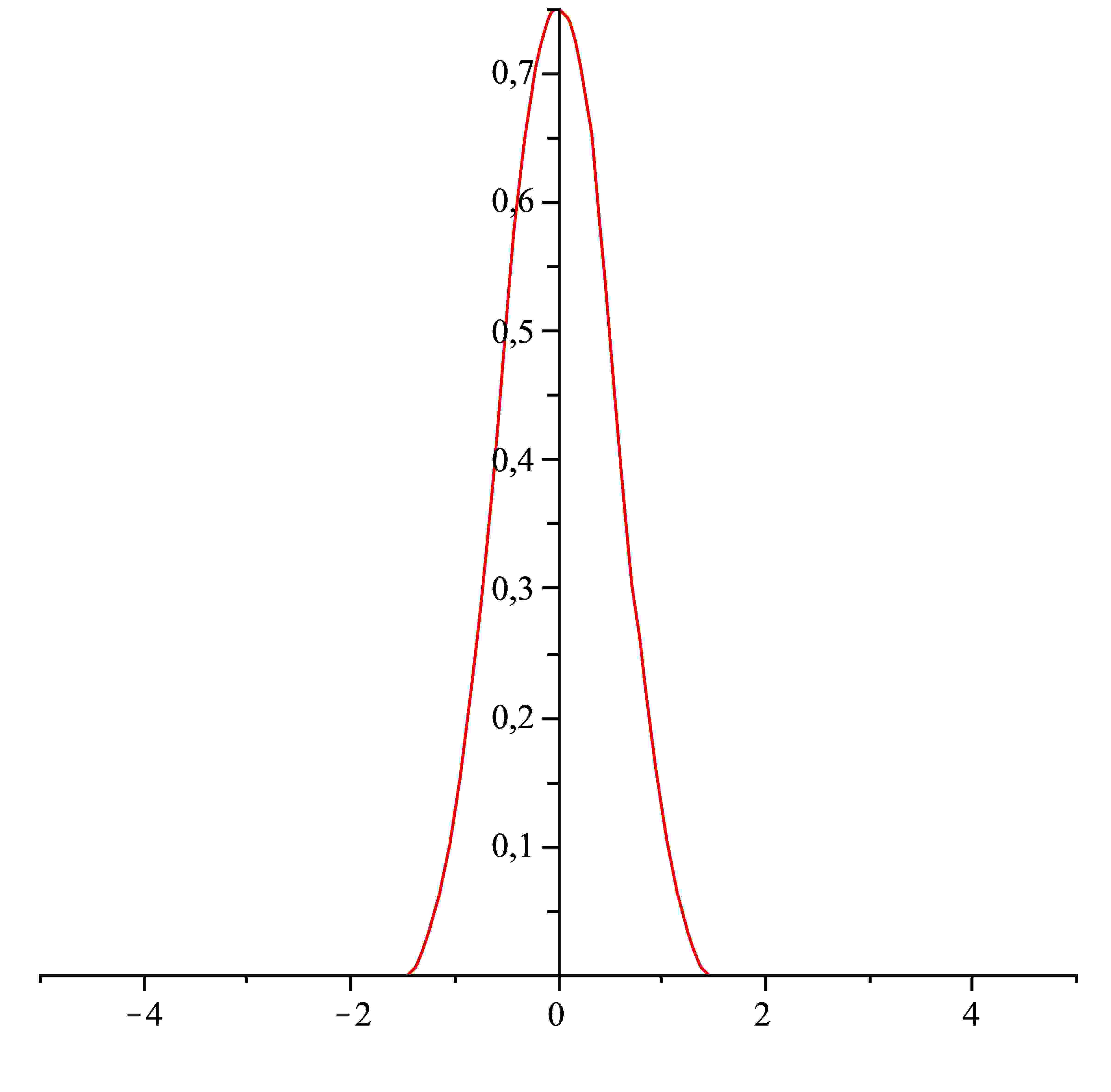}\includegraphics[width=.55\textwidth,height=.34\textwidth]
{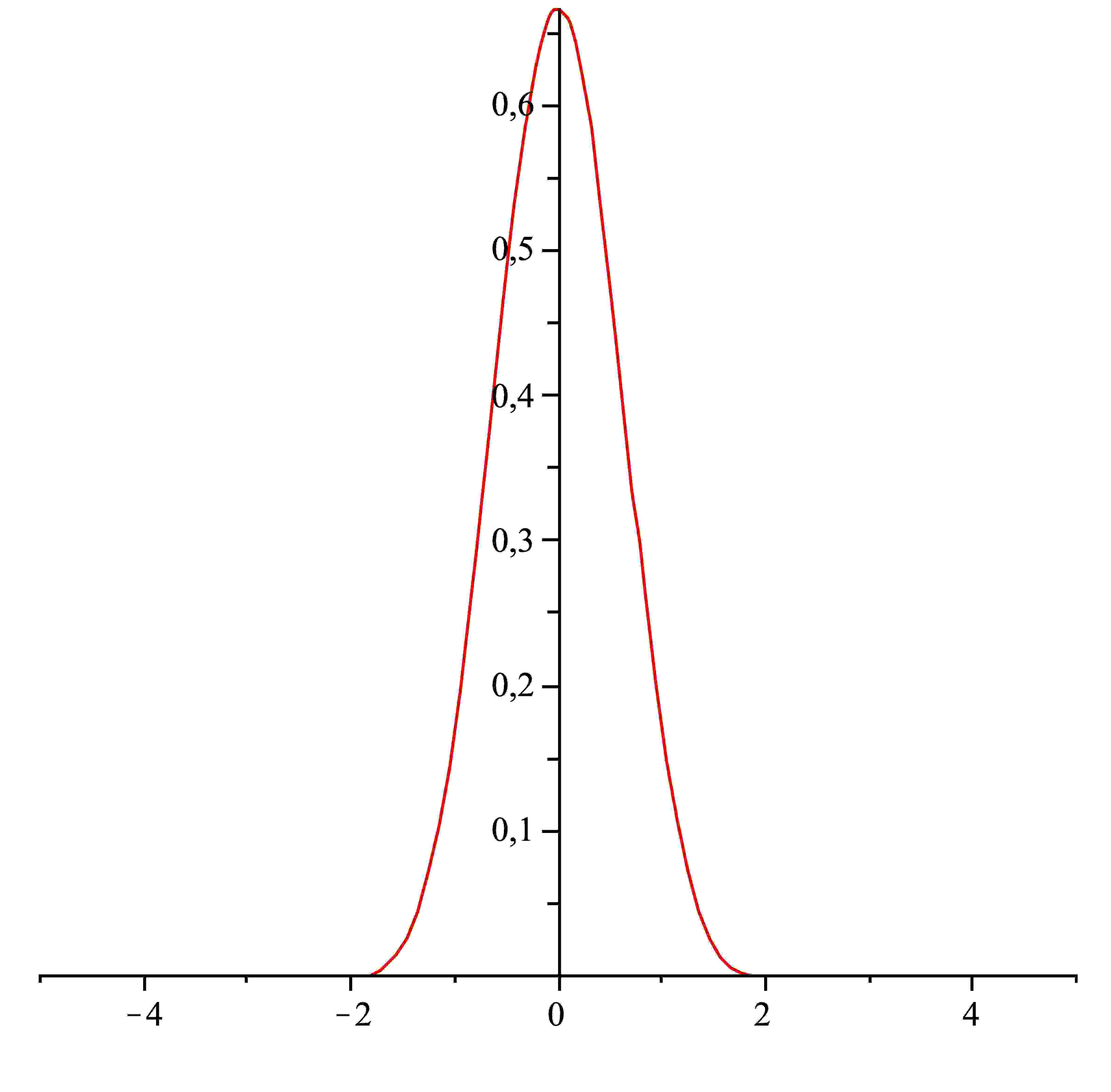}\end{center}
\vspace{-.5cm}
\caption{\small{\it The graphs of $M_3(u)$ and of $M_4(u)$ for $-5\leq u\leq5$.}}
\end{figure}
\vspace{-.7cm}
\begin{figure}[htbp]
\begin{center}\includegraphics[width=.7\textwidth,height=.47\textwidth]
{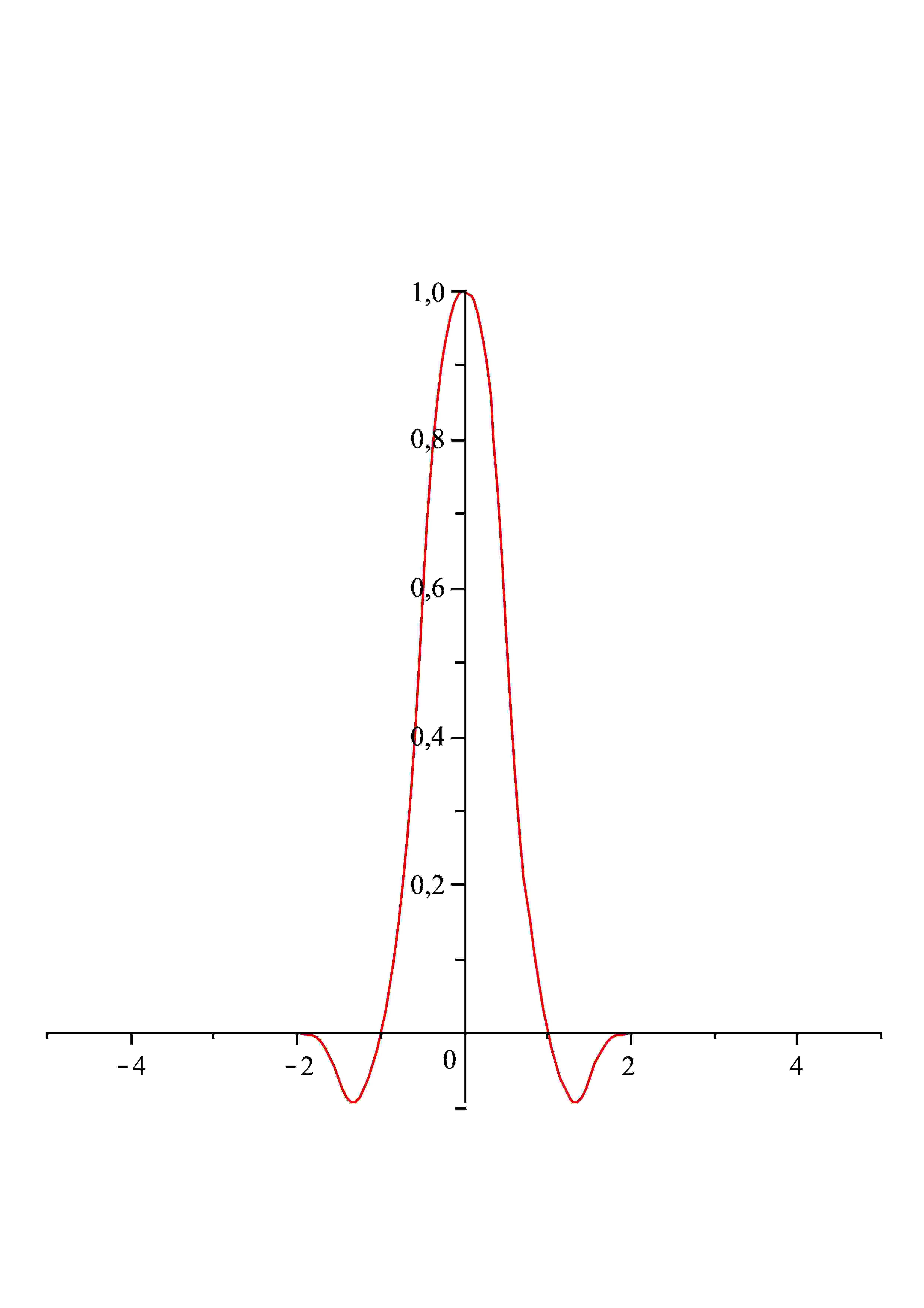}\end{center}

\vspace{-1.3cm}

\caption{\small{\it The graph of $M(u)$ for $-5\leq u\leq5$.}}
\end{figure}

In the following examples we will use the function $M(u)$ to define kernels which satisfy ($\chi_1$)--($\chi_6$).  We start from the example considered  in (2) of  Section 4. We choose the kernels $\chi_w(u)$ defined as $\chi_w(u)=M(wu).$ We have $$S_w^{(2)}f(x)=\int_{-\infty}^\infty M(wx-t)\left(\dfrac{w}{2}\int_{(t-1)/w}^{(t+1)/w}f(s)ds\right)dt.$$ 
We take a discontinuous function $f(x)$, as follows
$$f(x)=\begin{cases}3e^x,&x\mi-1,\\-1,&-1\leq x\mi0,\\2,&0\leq x\mi 1,\\ x,&1\leq x\mi 2,\\-2e^{-x},&x\geq 2.\end{cases}.$$ The graphs below (Figure 3) show the behaviour of the operator $S_w^{(2)}$ for $w=5,10$ and 15 respectively.
\begin{figure}[htbp]
\includegraphics[width=.49\textwidth,height=.49\textwidth]
{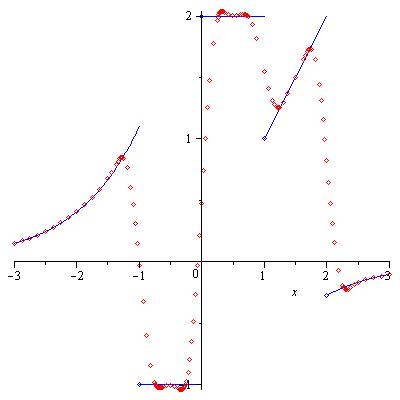}
\end{figure}

\begin{figure}[htbp]
\includegraphics[width=.49\textwidth,height=.49\textwidth]
{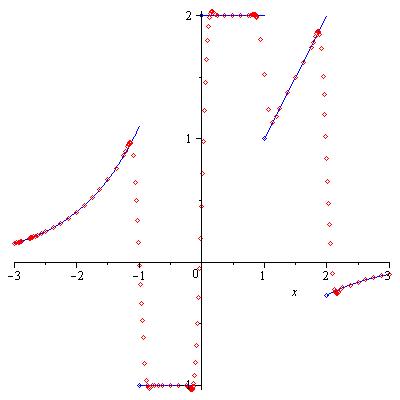}
\includegraphics[width=.49\textwidth,height=.49\textwidth]
{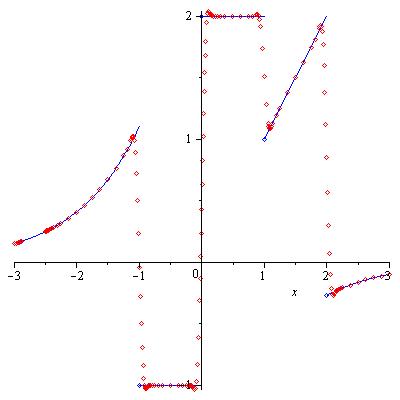}
\caption{\small{\it The graphs of the functions $S_5^{(2)}f(x),S_{10}^{(2)}f(x)$, $S_{15}^{(2)}f(x)$ compared with the graph of $f(x)$}}
\end{figure}

The next example represents the approximation of the operator $S_w^{(3)}$ of Section 4. In order to satisfy the assumptions ($\chi_1$)--($\chi_6$), we take the kernels $\chi_w(u)$ defined as $\chi_w(u)=wM(wu).$ So we are left with  $$S_w^{(3)}f(x)=\int_{\mathbb R}wM(wx-wt)\left(\dfrac{w}{2}\int_{t-\frac{1}{w}}^{t+\frac{1}{w}}f(s)ds\right)dt.$$ We take the function $$f(x)=\begin{cases}3e^x,&x\mi-1,\\-1,&-1\leq x\mi0,\\2,&0\leq x\mi 2,\\-2e^{-x},&x\geq 2.\end{cases}$$

\vspace{-.4cm}

Below we represent  the graphs of $S_5^{(3)}f(x), S_{10}^{(3)}f(x)$ and $S_{15}^{(3)}f(x)$ respectively (Fi\-gure 4).

\vspace{-.3cm}

\begin{figure}[htbp]
\includegraphics[width=.49\textwidth,height=.49\textwidth]
{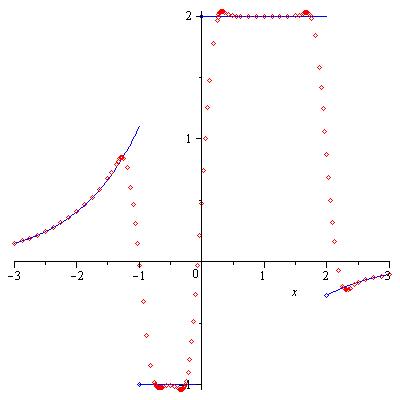}
\end{figure}

\vspace{-1cm}

\begin{figure}[htbp]
\includegraphics[width=.49\textwidth,height=.49\textwidth]
{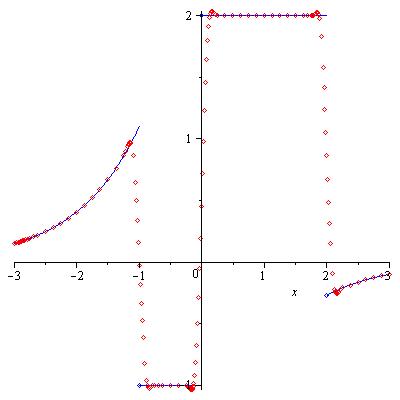}
\includegraphics[width=.49\textwidth,height=.49\textwidth]
{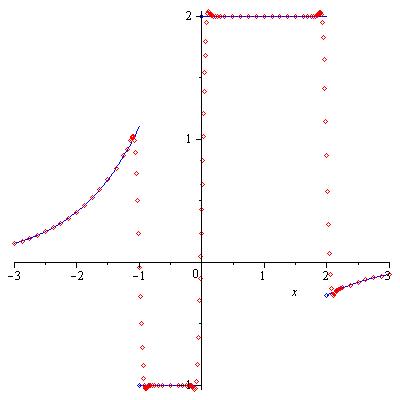}
\caption{\small{\it The graphs of the functions $S_5^{(3)}f(x),S_{10}^{(3)}f(x)$, $S_{15}^{(3)}f(x)$ compared with the graph of $f(x)$}}
\end{figure}

\vspace{3cm}

Our last graphical representation takes into account the example (4) in Section 4. In this case, however, we cannot take kernels based on the function $M(u)$ as before, because of the base space $\mathbb R^+$ and of the measure $d\mu(t)=\dfrac{dt}{t}$. Suitable kernel functions in this case are familiar to those working with Mellin operators, namely we consider the kernels 

\vspace{-.2cm}

$$\mathcal M_w(u)=\begin{cases} wu^w,&0<u<1,\\0,& \mbox{otherwise}\end{cases}.$$ 

\vspace{-.2cm}

Next, we consider the operators 
\vspace{-.2cm}
$$S_w^{(4)}f(x)=\int_0^\infty\mathcal M_w\left(\dfrac{x}{t}\right)\dfrac{1}{2\ln(1+1/w)}\left(\int_{t\frac{w}{w+1}}^{t\frac{w+1}{w}}f(u)\dfrac{du}{u}\right)\dfrac{dt}{t}.$$ 
\vspace{-.2cm}
We consider the function $$f(x)=\begin{cases}2x,&0\leq x<2,\\1,&2\leq x<4,\\-25/x^3,&x\geq 4\end{cases}.$$ 

\vspace{-.3cm}

Below (Figure 5) we represent the approximation of the operators  $S_w^{(4)}$ for $w=5,20$ and 30 respectively.
\begin{figure}[htbp]
\includegraphics[width=.49\textwidth,height=.49\textwidth]
{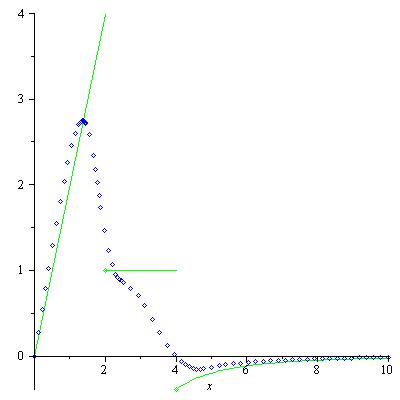}
\end{figure}

\vspace{-.8cm}

\begin{figure}[htbp]
\includegraphics[width=.49\textwidth,height=.49\textwidth]
{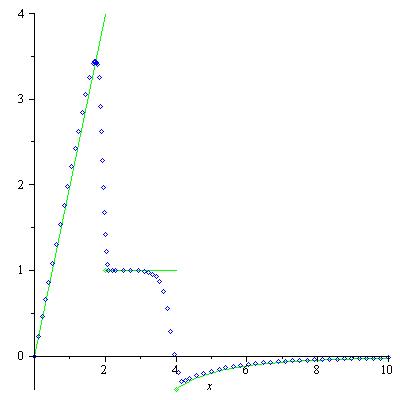}
\includegraphics[width=.49\textwidth,height=.49\textwidth]
{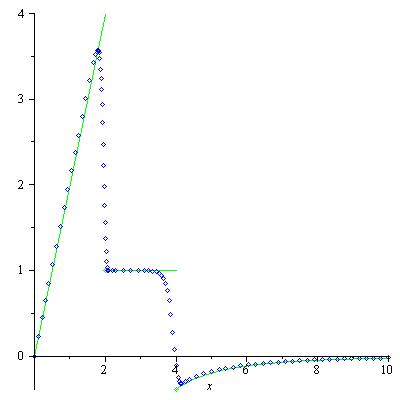}
\caption{\small{\it The graphs of the  functions $S_5^{(4)}f(x),S_{20}^{(4)}f(x)$, $S_{30}^{(4)}f(x)$ compared with the graph of $f(x)$}}
\end{figure}

\end{document}